\documentclass[MSNbibl,number,citesort,seceqn,dvips]{arxbj}
\usepackage{graphicx}
\aid{0}
\volume{19}
\issue{3}
\pubyear{2013}
\firstpage{931}
\lastpage{953}
\doi{10.3150/12-BEJ439} %kopijuoti is PTS

\makeatletter
\newcommand{\rrvert}{\vert}
\newcommand{\llvert}{\vert}
\newtheorem{lemma}{Lemma}[section]
\newtheorem{theorem}{Theorem}[section]
\newtheorem{corollary}{Corollary}[section]
\newcommand{\ess}{\operatorname{ess}}
\newcommand{\GARCH}{\operatorname{GARCH}}

\def\PROB{\mathbf{P}}
\def\EXP{\mathbf{E}}

\newcommand{\R}{\mathbb{R}}
\newcommand{\N}{\mathbb{N}}
\newcommand{\Z}{\mathbb{Z}}

\renewcommand{\P}{\mathcal{P}}
\newcommand{\B}{\mathcal{B}}
\newcommand{\F}{\mathcal{F}}
\newcommand{\T}{\mathcal{T}}
\renewcommand{\epsilon}{\varepsilon}
\makeatother

\begin{document}
\begin{frontmatter}

\title{On data-based optimal stopping under stationarity and ergodicity}
\runtitle{On data-based optimal stopping}

\begin{aug}
%%%% inicialai - be tarpu
\author[1]{\fnms{Michael} \snm{Kohler}\corref{}\thanksref{1}\ead[label=e1]{kohler@mathematik.tu-darmstadt.de}} \and
\author[2]{\fnms{Harro} \snm{Walk}\thanksref{2}\ead[label=e2]{walk@mathematik.uni-stuttgart.de}}
\runauthor{M. Kohler and H. Walk} %% auto
\address[1]{Fachbereich Mathematik, Technische Universit\"at Darmstadt,
Schlossgartenstr. 7, 64289 Darmstadt, Germany. \printead{e1}}
\address[2]{Fachbereich Mathematik, Universit\"at Stuttgart,
Pfaffenwaldring 57, 70569 Stuttgart, Germany.\\ \printead{e2}}
\end{aug}

% HISTORY:
\received{\smonth{9} \syear{2010}}
\revised{\smonth{4} \syear{2011}}

% ABSTRACT
%
\begin{abstract}
The problem of optimal stopping with finite horizon
in discrete time is considered in view of maximizing the expected
gain.
The algorithm proposed in this paper is completely
nonparametric in the sense that it
uses observed data from the past of the process up to time $-n+1$, $n
\in\N$,
not relying on any specific model assumption. Kernel regression estimation
of conditional expectations and prediction theory of individual
sequences are used as tools.
It is shown that the algorithm
is universally consistent:
the achieved expected gain converges to the
optimal value
for
$n \to\infty$
whenever the underlying process
is stationary and ergodic. An application to exercising American options
is given, and the algorithm is illustrated by simulated data.
\end{abstract}

% KEYWORDS
%

\end{frontmatter}

%s1 #&#
\section{Introduction}\label{se1}

In this paper an optimal stopping problem with finite horizon $L$
in discrete time is treated. The problem is formulated as follows:
Let
$(Z_j)_{j \in\Z}$
be a sequence of real-valued random variables and let
$g_j(Z_1^j)$,
with measurable bounded and real-valued functions $g_j$ on $\R^j$
$(g_0 = \mathit{const})$ and notation
$Z_1^j=(Z_1, \dots, Z_j)$,
be the gain when stopping at time
$j$ $(j \in\{0,1,\dots,L\})$.
In case that one stops at time $k \in\{0,1,\dots,L\}$
any stopping rule can rely only on the observed values
of $Z_j$ at times $j \leq k$.
Therefore, it can be described by a stopping time $\tau$, that is,
by a measurable function of
$Z_{-\infty}^L:=(\dots, Z_{-1}, Z_0, \dots, Z_L)$
where the event
$[\tau=k]$
is contained in the $\sigma$-algebra
$\F(Z_{-\infty}^k)$ generated by $Z_{-\infty}^k$.
Let
$\T(0,1,\dots,L)$ be the set of all such stopping times.
Any stopping time $\tau\in\T(0,1,\dots,L)$
yields the expected gain
\[
\EXP \bigl\{ g_\tau \bigl(Z_1^\tau \bigr) \bigr
\},
\]
and it is this quantity which one wants to maximize, that is,
one wants to construct a stopping time $\tau^* \in\T(0,1,\dots,L)$
such that
\[
V_0^* := \sup_{\tau\in\T(0,\dots,L)} \EXP \bigl\{ g_\tau
\bigl(Z_1^\tau \bigr) \bigr\} = \EXP \bigl\{
g_{\tau^*} \bigl(Z_1^{\tau^*} \bigr) \bigr\}
\]
(so-called value of the optimal stopping problem).
In the sequel, we assume only stationarity and ergodicity and define
decision rules on the basis of observed data. The unknown
underlying distribution is not used.

More precisely, for $n \in\N$ we assume that from the past of the process
the random variables $Z_{-n+1}, \dots, Z_0$ are observed and we want to
construct a stopping time
\[
\hat{\tau}_n= \hat{\tau}_n \bigl(Z_{-n+1}^L
\bigr) \in \T(0,1,\dots, L )
\]
such that
\[
\hat{V}_{0,n} := \EXP \bigl\{ g_{\hat{\tau}_n} \bigl(Z_1^{\hat{\tau}_n}
\bigr) \bigr\}
\]
converges to $V_0^*$.

In the definition of our estimates, we firstly use results
from the general theory of optimal stopping showing that
an optimal stopping time can be constructed
on the basis of dynamic programming
by recursively computing so-called
continuation value functions, which indicate the value of the
optimal stopping problem (from time $t$ on)
given an observed vector $Z_{-\infty}^t$
under the constraint of no stopping at time $t$
(cf., e.g., Chow, Robbins and Siegmund~\cite{RS71}
or Shiryayev~\cite{Sh78}). Secondly, we use that these continuation
values can be represented as conditional expectations
(cf., e.g., Tsitsiklis and van Roy~\cite{TsiRoy99}, Longstaff and Schwarz
\cite{LoSch99}
or Egloff~\cite{Eg05}), and our algorithm uses techniques from
nonparametric regression to estimate these conditional expectations
from observed
stationary and ergodic data
$Z_{-n+1}$, $Z_{-n}$, \ldots.
In contrast to the above references which study regression-based
Monte-Carlo methods for pricing American options,
for our estimates we do not use simulations of the underlying process,
because in our case its distribution is unknown, but use only
the observation of the individual sequence back to time $-n+1$.
This is in general a rather challenging
task, where usually extremely complex and data consuming algorithms
are necessary (cf., e.g., Morvai, Yakowitz and Gy\"orfi~\cite{MYG96}). But in
case that it is
enough to construct algorithms which converge in the so-called
Ces\`aro sense, a relatively simple and nice algorithm exists
(cf., e.g., Section 27.5 in Gy\"orfi \textit{et al.}~\cite{GyKoKrWa02}), which uses
techniques from the theory of prediction of individual
sequences (cf., e.g., Cesa-Bianchi and Lugosi~\cite{CBL06}).
These techniques have already been
used successfully in the context of portfolio
optimization (cf., e.g.,
Gy\"orfi, Lugosi and Udina~\cite{GLU06},
Gy\"orfi, Udina and Walk~\cite{GUW08} and the references therein).
In this paper, we introduce as main
trick an averaging of such estimates
and show that by using this trick we can derive
a consistency result of our estimated stopping rule from Ces\`aro
consistency of the underlying regression estimates.
So in the definition of our estimate, we thirdly apply estimates
defined by use of ideas from the prediction theory of individual
sequences.

As an application, we consider the problem of exercising an American option
in discrete time (also called Bermudan option) in view of maximizing
of the expected discounted payoff.

The algorithm computing
estimates of the optimal
stopping time is described
in Section~\ref{se2} and the main result is formulated in Section~\ref{se3}, where also
an application to American options is described. In Section~\ref{se4}, we
illustrate our algorithm by applying it to simulated data,
Section~\ref{se5} contains the proof of the main result,
the proof of an auxiliary result is given in the \hyperref[app]{Appendix}.
%in Section 2, the main result is formulated in Section 3, and
%the algorithm is illustrated by applying it to simulated
%data in Section 4. Section 5 contains the proof of the main result,
%proofs of auxiliary results are given in the appendix.

%s2 #&#
\section{Construction of an approximation of the optimal stopping~time}\label{se2}

Our first idea is to use results from the general theory of optimal
stopping in order to determine the optimal stopping time $\tau^*$.
Let $t \in\{0,\dots,L-1\}$ be fixed and denote the set of all
stopping times
with values in $\{t+1,\dots,L\}$ by $\T(t+1,\dots,L)$.

For each $\tau\in\T(t+1,\dots,L)$, define the real random variable
$h_\tau$ on the probability space
$(\prod_{-\infty}^t \R, \bigotimes_{-\infty}^t \B, \PROB_{Z_{-\infty}^t})$
with Borel product $\sigma$-algebra and distribution of $Z_{-\infty}^t$
by
\[
h_\tau \bigl(z_{-\infty}^t \bigr) := \EXP \bigl\{
g_\tau \bigl(Z_1^\tau \bigr) |
Z_{-\infty}^t = z_{-\infty}^t \bigr\}.
\]
Then, according to Chow, Robbins and Siegmund~\cite{RS71}, Section 7.6,
\[
\ess \sup_{\tau\in\T(t+1,\dots,L)} h_\tau=: q_t
\]
is defined as a real-valued random variable $y$ on this probability
space such that
\begin{enumerate}[(ii)]
\item[(i)]
$\PROB_{Z_{-\infty}^t} \{ y \geq h_\tau\} =1$
for every $\tau\in\T(t+1,\dots,L)$,

\item[(ii)]
if $y^\prime$ is any real random variable on the probability space
satisfying
\[
\PROB_{Z_{-\infty}^t} \bigl\{ y^\prime\geq h_\tau \bigr\} =1
\qquad \mbox{for every } \tau\in\T(t+1,\dots,L),
\]
then
\[
\PROB_{Z_{-\infty}^t} \bigl\{ y^\prime\geq y \bigr\} =1.
\]
\end{enumerate}
Thus, $q_t$ is unique $mod$ $\PROB_{Z_{-\infty}^t}$, that is,
two versions of $q_t$ coincide $\PROB_{Z_{-\infty}^t}$-almost everywhere.
By Theorem 1.5 in Chow, Robbins and Siegmund~\cite{RS71}, $q_t$ always
exists, and there exists a countable subset $\T_t^*$ of $\T(t+1,\dots,L)$
such that
\[
q_t = \sup_{\tau\in\T_t^*} h_\tau.
\]
Furthermore we set $q_L:=0$. $q_t$ is denoted as continuation value
function $(t \in\{ 0, \dots, L\})$.

The so-called continuation values
\begin{eqnarray*}
q_t \bigl(z_{-\infty}^t \bigr)&\hspace*{2.8pt}=:&
\ess \sup_{\tau\in\T(t+1, \dots, L)} \EXP \bigl\{ g_\tau \bigl(Z_1^\tau
\bigr) | Z_{-\infty}^t = z_{-\infty}^t \bigr\}
\qquad \bigl(t \in\{0, \dots, L-1\} \bigr),
\\
q_L \bigl(z_{-\infty}^L \bigr)&=&0
\end{eqnarray*}
describe the values of the optimal stopping problem
from $t$ on given $Z_{-\infty}^t = z_{-\infty}^t$
subject to the constraint of not stopping at time $t$.

Replacing $\T(t+1, \dots, L)$ by $\T(t, \dots, L)$ leads
to the so-called value functions
%
%e2.1 #&#
\begin{equation}
\label{weq1} \ess \sup_{\tau\in\T(t, \dots, L)} h_\tau=: V_t\qquad
\bigl(t \in\{ 0, \dots, L\} \bigr).
\end{equation}
$V_t(z_{-\infty}^t)$ describes the value of the optimal
stopping problem (from $t$ on) given $Z_{-\infty}^t = z_{-\infty}^t$.

For $t \in\{-1,0, \dots, L-1\}$, set
%
%e2.2 #&#
\begin{equation}
\label{weq2} \tau_t^* := \inf \bigl\{ s \geq t+1 \dvt
q_s \bigl(Z_{-\infty}^s \bigr) \leq g_s
\bigl(Z_1^s \bigr) \bigr\}.
\end{equation}

We can conclude from the general theory
of optimal stopping (see, e.g., Chow, Robbins and Siegmund~\cite{RS71}
or Shiryayev~\cite{Sh78}).

%le2.1 #&#
\begin{lemma}
\label{le1}
It holds
%
%e2.3 #&#
\begin{equation}
\label{le1eq1} V_t \bigl(z_{-\infty}^t \bigr) =
\EXP \bigl\{ g_{\tau_{t-1}^*} \bigl(Z_1^{\tau_{t-1}^*} \bigr) |
Z_{-\infty}^t=z_{-\infty}^t \bigr\}
\end{equation}
$\PROB_{Z_{-\infty}^t}$-almost everywhere
for $t \in\{0, \dots, L\}$.
Furthermore
%
%e2.4 #&#
\begin{equation}
\label{le1eq2} V_0^* := \sup_{\tau\in\T(0, \dots, L)} \EXP \bigl\{
g_\tau \bigl(Z_1^\tau \bigr) \bigr\} = \EXP
\bigl\{ g_{\tau^*} \bigl(Z_1^{\tau^*} \bigr) \bigr\}
\end{equation}
is fulfilled for
\[
\tau^* := \tau^*_{-1} = \inf \bigl\{ j \in\{ 0, 1, \dots, L\} \dvt
g_j \bigl(Z_1^j \bigr) \geq q_{j}
\bigl(Z_{-\infty}^j \bigr) \bigr\}.
\]
\end{lemma}
Lemma~\ref{le1} can be proven as in the case of
Markovian processes (cf., e.g., proof of Theorem 1 in Kohler~\cite{Ko10}),
a complete proof of this lemma is available from the authors by request.

From Lemma~\ref{le1}, we get that it suffices to compute the
continuation value functions $q_{0}$, \dots, $q_{{L-1}}$ in order
to construct the optimal stopping rule $\tau^*$.
In Tsitsiklis and van Roy~\cite{TsiRoy99}, Longstaff and Schwarz
\cite{LoSch99}
and Egloff~\cite{Eg05} it is shown that in case of Markovian
processes the continuation values can be computed recursively
by evaluation of conditional expectations.
The same can be shown also in the setting considered
in this paper.

%le2.2 #&#
\begin{lemma}\label{le2}
The continuation values satisfy
%
%e2.5 #&#
\begin{equation}
\label{le2eq2} q_j \bigl(z_{-\infty}^j \bigr) = \EXP
\bigl\{ \max \bigl\{ g_{j+1} \bigl(Z_1^{j+1}
\bigr),q_{j+1} \bigl(Z_{-\infty}^{j+1} \bigr) \bigr\} |
Z_{-\infty}^j = z_{-\infty}^j \bigr\}
%&=&
%g_{j+1}(z_1^{j},Z_{j+1}),q_{j+1}(z_{-\infty}^{j},Z_{j+1})
%Z_{-\infty}^j = z_{-\infty}^j
\end{equation}
$\PROB_{Z_{-\infty}^j}$-almost everywhere
and
%
%e2.6 #&#
\begin{equation}
\label{le2eq3} q_j \bigl(z_{-\infty}^j \bigr) =
\EXP \bigl\{ g_{\tau_j^*} \bigl(Z_1^{\tau_{j}^*} \bigr) |
Z_{-\infty}^j = z_{-\infty}^j \bigr\}
\end{equation}
$\PROB_{Z_{-\infty}^j}$-almost everywhere
for any $j \in\{0,1, \dots, L-1\}$.
\end{lemma}
Lemma~\ref{le2} can be proven as in the case of
Markovian processes (cf., e.g., proof of Theorem 2 in Kohler~\cite{Ko10}),
again a complete proof of this lemma is available
from the authors by request.

Usually in applications, the distribution of the
underlying process $(Z_n)_n$
is unknown and therefore it is impossible to use
(\ref{le2eq2}) (or (\ref{le2eq3}))
in order to compute the continuation values.
In the sequel, we will try to estimate them by using
(recursively defined) regression estimates in order to
approximate the conditional expectations in (\ref{le2eq2}).
To do this, for any $n \in\N$ we use
$Z_{-n+1}^0$
in order to
construct an estimate of the optimal stopping rule
on the data $Z_0$, \dots, $Z_L$.

Next, we describe how we construct estimates
$
\hat{q}_j^{(n)}( Z_{-n+1}^j)
$
of $q_j(Z_{-\infty}^j)$.

The estimates are defined recursively with respect
to $j \in\{0, \dots, L\}$. For $j=L$, we have $q_L=0$ and in
this case we set
\[
\hat{q}_L^{(n)}:=0.
\]
Given $\hat{q}_{j+1}^{(m)}$
(defined on $\R^{j+m+1}$), $m \leq n$,
for some $j \in\{0,1, \dots, L-1\}$ we
define $\hat{q}_{j}^{(n)}$
as follows.

To make the construction more transparent, for the function
$q_j\dvtx \prod_{-\infty}^j \R\rightarrow\R$, which by (\ref{le2eq2})
is given as a regression function, we define a regression estimation
function $\hat{m}^{(n)}_{j,(k,h)}(z_{-n+1}^0; \cdot)\dvtx \prod_{-n+1}^j \R
\rightarrow
\R_+$
with parameters $k$, $h$ using realizations $z_{-n+1}^0$ of $Z_{-n+1}^0$.
The definition
depends on parameters
$k \in\N$ (indicating how far back the estimate will look, and thus
indicating also the dimension of the occurring regression estimation problem)
and $h>0$ (a so-called bandwidth
which (roughly speaking) indicates how similar observed values
in the past must be to the current observed values in order to be
included in the prediction of the future value)
and a kernel function $K\dvtx \R^{j+k+1} \rightarrow\R_+$.
We define the latter by
\[
K(v) := H \bigl( \|v\|_2^{j+k+1} \bigr),
\]
where $\|v\|_2$ denotes the Euclidean norm of $v$ and $H\dvtx \R_+
\rightarrow\R_+$
is a given nonincreasing and continuous function satisfying
\[
H(0)>0 \quad\mbox{and} \quad t \cdot H(t) \rightarrow0 \qquad(t \rightarrow
\infty)
\]
(e.g., $H(v) = \mathrm{e}^{-v^2}$).
The use of the exponent $j+k+1$ in the definition of $K$ and not for the
factor $t$ in the condition on $H$ allows to choose $H$ independent
of $j$ and $k$.

We set $\hat{m}_{L,(k,h)}^{(n)}(z_{-n+1}^0; \cdot):=0$ and
use local averaging to define
%
%e2.7 #&#
\begin{eqnarray}
\label{eqw*} && \hat{m}_{j,(k,h)}^{(n)} \bigl(z_{-n+1}^0;u_{-n+1}^j
\bigr)
\nonumber
\\
&&\quad := \sum_{i=-n+k+1}^{-(j+1)} \max \bigl\{
g_{j+1} \bigl(z_{i+1}^{i+j+1} \bigr),
\hat{q}_{j+1}^{(n+i)} \bigl(z_{-n+1}^{i+j+1} \bigr)
\bigr\}
\\
&& \hspace*{63pt}{}\cdot { K \biggl( \frac{u_{-k}^j-z_{i-k}^{i+j}}{h} \biggr) }\bigg /{ \sum
_{l=-n+k+1}^{-(j+1)} K \biggl( \frac{u_{-k}^j-z_{l-k}^{l+j}}{h} \biggr) }
\nonumber
\end{eqnarray}
for
$u_{-n+1}^j \in\prod_{-n+1}^j \R$, where
we set
\[
\hat{m}_{j,(k,h)}^{(n)} \bigl(z_{-n+1}^0;
\cdot \bigr) :=0
\]
for $k \geq n-j-1$, and $\frac{0}{0}:=0$.
Then we set
\[
\hat{q}_{j,(k,h)}^{(n)} \bigl(z_{-n+1}^j
\bigr) := \hat{m}_{j,(k,h)}^{(n)} \bigl(z_{-n+1}^0;z_{-n+1}^j
\bigr).
\]

Let $h_r>0$ be such that $h_r \rightarrow0$ for $r \rightarrow\infty$
and set
\[
\P:= \bigl\{ (k,h_r) \dvt k,r \in\N \bigr\}.
\]
For $(k,h) \in\P$ define the cumulative loss of the corresponding
estimate by
%
%e2.8 #&#
\begin{eqnarray}
\label{eqw**} \hat{Q}_{n,j}(k,h) &:=&\hat{Q}_{n,j}
\bigl(z_{-n+1}^j,k,h \bigr)
\nonumber
\\
&:=& \frac{1}{n} \sum_{i=1}^{n-1}
\bigl( \hat{q}_{j,(k,h)}^{(i)} \bigl( z_{-n+1}^{-n+i+j}
\bigr)
\\
&&\hspace*{26pt} {} - \max \bigl\{ g_{j+1} \bigl( z_{-n+i+1}^{-n+i+j+1}
\bigr), \hat{q}_{j+1}^{(i)} \bigl( z_{-n+1}^{-n+i+j+1}
\bigr) \bigr\} \bigr)^2.
\nonumber
\end{eqnarray}

Put $c = 8 B^2$ (where we assume that the gain functions
are bounded by $B$), let $(p_{k,r})_{k,r}$
be a probability distribution such that $p_{k,r}>0$
for all $k,r \in\N$, and define weights, which depend
on these cumulative losses, by
\[
w_{n,k,r}^{(j)} := w_{n,k,r}^{(j)}
\bigl(z_{-n+1}^j \bigr) :=p_{k,r} \cdot
\mathrm{e}^{-n \cdot\hat{Q}_{n,j}(k,h_r)/c}
\]
and their normalized values
\[
v_{n,k,r}^{(j)}:= v_{n,k,r}^{(j)}
\bigl(z_{-n+1}^j \bigr) := \frac{w_{n,k,r}^{(j)}}{
\sum_{s,t=1}^\infty
w_{n,s,t}^{(j)}
}.
\]
The estimate $
\hat{q}_j^{(n)}
$
is defined on $\prod_{-n+1}^j \R$
as the convex combination of the estimates
$\hat{q}_{j,(k,h_r)}^{(n)}$
using the weights
$v_{n,k,r}^{(j)}$, that is, $
\hat{q}_j^{(n)}
$
is defined
by
%
%e2.9 #&#
\begin{equation}
\label{eqw***} \hat{q}_j^{(n)} \bigl(z_{-n+1}^j
\bigr) := \sum_{k,r=1}^\infty
v_{n,k,r}^{(j)} \cdot \hat{q}_{j,(k,h_r)}^{(n)}
\bigl(z_{-n+1}^j \bigr).
\end{equation}
Finally, for the computation of our estimated stopping rule
we use the arithmetic mean of the first $n$ estimates,
that is, we use
%
%e2.10 #&#
\begin{equation}
\label{se2eqaveraging} \hat{q}_{j,n} \bigl(z_{-n+1}^j
\bigr) := \frac{1}{n} \sum_{l=1}^n
\hat{q}_j^{(l)} \bigl( z_{-l+1}^j \bigr)
\end{equation}
for $j \in\{0,1, \dots, L-1\}$ and
$\hat{q}_{L,n}:=q_{L}=0$.\vadjust{\goodbreak}

With this estimate of $q_j$, we estimate the optimal stopping rule
\[
\tau^* :=\inf \bigl\{ j \in\{0, 1, \dots, L\} \dvt g_j
\bigl(Z_1^j \bigr) \geq q_{j}
\bigl(Z_{-\infty}^j \bigr) \bigr\}
\]
by
\[
\hat{\tau}_n :=\inf \bigl\{ j \in\{0, 1, \dots, L \} \dvt
g_j \bigl(Z_1^j \bigr) \geq
\hat{q}_{j,n} \bigl(Z_{-n+1}^j \bigr) \bigr\}.
\]

%s3 #&#
\section{Main theoretical result}\label{se3}

In Theorem~\ref{th1} below, we assume
that the underlying process
$(Z_j)_{j \in\Z}$ in $\R$ is (strictly) stationary and ergodic,
that is, for each $B \in\B_\Z$ (where $\B_\Z$ is the Borel $\sigma
$-algebra
in $\R^\Z$) and each $k \in\Z$
\[
\PROB \bigl\{ (Z_j)_{j \in\Z} \in B \bigr\} = \PROB \bigl\{
(Z_{j+k})_{j \in\Z} \in B \bigr\}
\]
and for each $B \in\B_\Z$ such that the event
\[
A:= \bigl\{ (Z_{j+k})_{j \in\Z} \in B \bigr\}
\]
does not depend on $k \in\Z$ one has
\[
\PROB(A) \in\{0,1\}
\]
(cf., e.g., G\"anssler and Stute~\cite{GS77} or
Gy\"orfi \textit{et al.}~\cite{GyKoKrWa02}, page 565).

Let the estimate $\hat{\tau}_n$ of the optimal stopping rule
$\tau^*$ be defined as in the previous section. Then the
following result is valid.

%th3.1 #&#
\begin{theorem}
\label{th1}
Let $(Z_j)_{j \in\Z}$ be an arbitrary stationary and ergodic
sequence of real-valued
random variables. Assume
that the gain functions $g_l\dvtx \R^l \rightarrow\R$ $(l=0, \dots, L)$
with $g_0 = \mathit{const}$
are measurable,
nonnegative and bounded (in absolute value) by $B>0$.
Let the estimate be defined as in Section~\ref{se2}, where the
kernel $K$ is given by
\[
K(v) = H \bigl( \|v\|_2^{j+k+1} \bigr)
\]
for some $H\dvtx \R_+ \rightarrow\R_+$ which
is a nonincreasing and continuous function satisfying
\[
H(0)>0 \quad\mbox{and} \quad t \cdot H(t) \rightarrow0 \qquad(t \rightarrow
\infty).
\]
Then
\[
\hat{V}_{0,n} := \EXP \bigl\{ g_{\hat{\tau}_n} \bigl(Z_1^{\hat{\tau}_n}
\bigr) \bigr\} \rightarrow V_0^* = \EXP \bigl\{ g_{\hat{\tau}^*}
\bigl(Z_1^{\hat{\tau}^*} \bigr) \bigr\}
\]
for $n \rightarrow\infty$.
\end{theorem}

As an application, we consider the problem of exercising an American option
in discrete time in view of maximization of the expected payoff. Let $X_j$,
$j \in\Z$, be positive random variables
defined on the same probability\vadjust{\goodbreak} space describing the values of
the underlying asset of the option at time points $j \in\Z$.
For simplicity, we consider
only the case that
$X_{j}$ be real-valued, that is, we consider only options on a single
asset.
Hereby,
we assume only that the corresponding returns
$Z_j:=X_j/X_{j-1}$
form a stationary
and ergodic sequence. The unknown underlying distribution
is not used.
Let $f\dvtx \R\rightarrow\R_+$ be the payoff function
of the option, which we assume to be nonnegative, bounded and measurable,
for example, $f(x)=\max\{K-x,0\}$
in case of an American put option with strike $K$. Let $r^*$ be the
riskless interest rate. If we get the payoff at time $t >0$,
we discount it towards zero by the factor $\mathrm{e}^{-r^* \cdot t}$, so
for asset value $x$ at time $t$
the discounted payoff of the
option is
$\mathrm{e}^{-r^* \cdot t} \cdot f(x)$.

Let $L>0$ be the expiration date of our option.
In the sequel, we renormalize the payoff function such that we can assume
$X_{0}=100$, and we consider
an American option on $X_j$ with exercise opportunities restricted
to $\{0,1,\dots,L\}$ (sometimes also called Bermudan option).
Any rule for exercising such an option
within $\{0,1,\dots,L\}$
can be described by a stopping time $\tau\in\T(0, \dots, L)$.
Any stopping time $\tau$
describing the exercising of an
American option yields in the mean the payoff
\[
\EXP \bigl( \mathrm{e}^{- r^* \cdot\tau} \cdot f(X_\tau) \bigr),
\]
which we want to maximize, that is, we want
to construct a stopping time $\tau^* \in\T({0}, \dots, {L})$
such that
\[
V_0^* := \sup_{\tau\in\T(0, \dots, L)} \EXP \bigl\{ \mathrm{e}^{- r^* \cdot\tau}
\cdot f(X_\tau) \bigr\} = \EXP \bigl\{ \mathrm{e}^{- r^* \cdot\tau^*} \cdot
f(X_{\tau^*}) \bigr\}.
\]
It should be noted that $V_0$ is not the price of the option as
defined in financial mathematics since we ignore the rest of
the financial market, in particular we do not buy, sell or borrow
additional stocks in parallel. Instead, we are dealing
with the situation of a holder of the option who has
no other possibilities than to exercise the option.

We assume that $X_{-n}, \dots, X_0$ or -- equivalently -- $Z_{-n+1}$, \ldots, $Z_0$ are observed. Then we set $g_0=f(X_0)=f(100)$,
$g_j(Z_1^j)=\mathrm{e}^{-r^* \cdot j} f(X_j)
=\mathrm{e}^{-r^* \cdot j} f(100 \cdot Z_1 \cdot Z_2 \cdots Z_j)$
$(j=1, \dots, L)$ and define the
sequence of stopping times $\hat{\tau}_n$ as in Section~\ref{se2}. Immediately
from Theorem~\ref{th1}, we can conclude the following corollary.

%co3.1 #&#
\begin{corollary}
Let $(X_j)_{j \in\Z}$ be an arbitrary sequence of positive
random variables such that the corresponding
returns are stationary and ergodic. Assume
that the payoff function is measurable,
nonnegative and bounded by $B>0$.
Let the estimate be defined as above, where the
kernel $K$ is given by
\[
K(v) = H \bigl( \|v\|_2^{j+k+1} \bigr)
\]
for some $H\dvtx \R_+ \rightarrow\R_+$ which
is a nonincreasing and continuous function satisfying
\[
H(0)>0 \quad\mbox{and} \quad t \cdot H(t) \rightarrow0 \qquad(t \rightarrow
\infty).
\]
Then
\[
\hat{V}_{0,n} := \EXP \bigl\{ \mathrm{e}^{- r^* \cdot\hat{\tau}_n} \cdot
f(X_{\hat{\tau}_n}) \bigr\} \rightarrow V_0^* = \EXP \bigl\{
\mathrm{e}^{- r^* \cdot\tau^*} \cdot f_{\tau^*}(X_{\tau^*}) \bigr\}
\]
for $n \rightarrow\infty$.\vadjust{\goodbreak}
\end{corollary}

%s4 #&#
\section{Application to simulated data}\label{se4}

\begin{figure}

\includegraphics{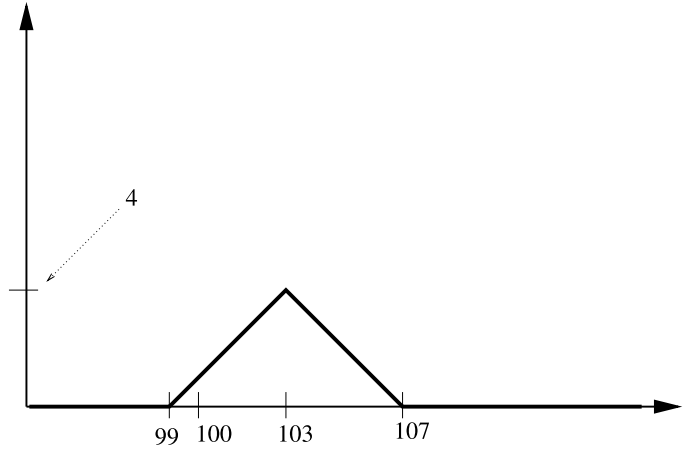}

\caption{Butterfly spread payoff function used in the
simulation.}\label{fig1}
\end{figure}

In this section, we evaluate the behaviour of our newly
proposed estimate for finite sample size
by applying it to simulated data. Here, we consider the optimal
exercising of an American option in discrete time
which can be exercised on one of the five equidistant time points
$t_0=0$, $t_1=0.25$, $t_2=0.5$, $t_3=0.75$
and $t_4=1$. The starting value of the stock is $x_0=100$,
for the payoff function we use a butterfly payoff function given by
$f(x)=\max\{0,\min\{x-99,107-x\}\}$ (cf., Figure~\ref{fig1}).
%
%
%f1 #&#
As model for generating the stock values, we consider a $\GARCH(1,1)$ model
in the form of Duan~\cite{Du95}. Here, we simulate the price process
according to
\begin{eqnarray*}
X_{i+1} &=& X_i \cdot\exp \biggl( \frac{r^*}{4} -
\frac{1}{2} \cdot\sigma_{i+1}^2 + \sigma_{i+1}
\cdot\epsilon_{i+1} \biggr),
\\
\sigma_{i+1}^2 &=& \delta_0 +
\delta_1 \cdot(\sigma_{i} \cdot\epsilon_{i} -
\lambda\cdot\sigma_i)^2 + \xi_1 \cdot
\sigma_i^2,
\end{eqnarray*}
where
$r^*=0.05$, $\lambda=0.7136$, $\delta_0=0.0000664$,
$\delta_1=0.144$, $\xi_1=0.776$ and where $(\epsilon_t)_{t \in\Z}$
are independent normally distributed random variables
with expectation zero and variance one. We start our simulation
with $X_0=x_0=100$. For $\sigma_0$, we use the random value we
get if we start the second recursion with $\sigma^2_{-1600}=0$.

We consider four different algorithms to estimate the optimal stopping
time: The first two algorithms are
simple methods where we exercise the option at the first time
when the payoff is greater than zero (\textit{simple1}) or at the
expiration date
of the option (\textit{simple2}). The third algorithm is the newly proposed
algorithm of this article (\textit{new algorithm}), and the fourth algorithm
(\textit{optstop})
is a regression-based Monte Carlo estimate of the optimal stopping rule
based on the true price process,
where we extend the state space in order to get a 3-dimensional
Markovian process (i.e., we use $(X_i, \sigma_i, \epsilon_i)$
as variables of the algorithm).
As regression-based Monte Carlo procedure, we use
the smoothing spline algorithm described in Kohler~\cite{Ko06}, which gives
results which are usually at least comparable but often better than the
algorithms of Tsitsiklis and Van Roy~\cite{TsiRoy99} and Longstaff and Schwarz
\cite{LoSch99}
based on parametric regression (cf. Kohler~\cite{Ko10}). This algorithm can
never be used in a real application since it
requires that the distribution of the underlying data is known
and since its decisions depend on the not
observable random variables $\sigma_i$ and $\epsilon_i$, however, it can
be considered as an approximation of the theoretical optimal stopping rule.

In contrast to algorithms one and two, the
algorithms three and four require training data. Our newly proposed
algorithm three uses a path of values of length $1500$
(which is part of the path of length $1600$ preceding our evaluation
paths)
in order to learn its stopping rule, that is, it depends on observable
values of the stock from the past. The theoretical algorithm four requires
a training set consisting of paths generated independently and
identically to
the path for which it should generate the stopping rule (which is
never available in any real application).
In our simulation, the algorithm is
based on $1000$ paths of length $5$ starting with $x_0=100$, each of
them extending the same path before time $t_0=0$ used also in the
evaluation of the stopping rule.

All algorithms are evaluated by applying them to $1000$ paths
of length $5$ starting with $x_0=100$, where each of
them extends the same path before time $t_0=0$, and we compute
the average of the $1000$ payoffs achieved. Since for two of our
four algorithms this result depends on the random training data,
we repeat this whole procedure $100$ times and report the means and the
standard deviations of the resulting values for each algorithm.

In the practical implementation of our newly proposed algorithm,
we consider as bandwidths $h \in\{0.001,0.01,0.1\}$ and use
the $k \in\{0,1,2\}$ last values of the returns for prediction of
the value at the next time step. Each of these $3 \cdot3 = 9$ models
gets the same probability $p_{k,r}=\frac{1}{9}$, and for the
constant used for computing the weights of the estimate from the
cumulative empirical losses we use $c=8 B^2$ where $B$ is the maximal
value of the payoff function.
In addition, we make the following modifications: Firstly, we simplify the
computation of the algorithm in such a way that we do not
use the final averaging step
(\ref{se2eqaveraging}), because otherwise we are not able
to compute the result of our algorithm in a reasonable time on
a standard computer.
Secondly, we do not use returns relative to the previous
day as $x$-values for our regression estimates, instead we
use returns relative to the beginning of the time interval
of an option. With the later modification, it can be shown that the
theoretical result above is still valid because
a consecutive sequence of these modified returns generates
the same $\sigma$-algebra as the corresponding original returns.
Finally, we ignore the first $n_{0,t}=(L-1-t) \cdot200$ data points
during the computation of $
\hat{q}_t^{(n)}$ since we think
that the first $n_{0,t}$ values of $\hat{q}_{t+1}^{(n)}$ are not reliable
because they are based on too few data points.

The results of the four algorithms are reported in Table~\ref{tab1}.
As we can see from Table~\ref{tab1} both simple algorithms are
clearly outperformed by our newly proposed algorithm, which achieves
results which are very close to the results of the
exercising strategy \textit{optstop}
relying on information not available in a real application.

%
%t1 #&#
\begin{table}
\tabcolsep=0pt
\caption{Achieved payoffs by the four different algorithms}\label{tab1}
\begin{tabular*}{\textwidth}{@{\extracolsep{\fill}}lllll@{}}
\hline
& Simple1&  Simple2&  New algorithm&  Optstop\\
\hline
Mean value & $\hphantom{(} 1.00$ & $\hphantom{(}0.61$ & $\hphantom{(}1.64$ & $\hphantom{(}1.72$ \\
(Standard deviation) & $(0.00)$& $(0.04)$ & $(0.47)$ & $(0.04)$\\
\hline
\end{tabular*}
\end{table}

From Table~\ref{tab1}, we see that in principle our new algorithm
could also be used as a numerical tool to evaluate American options
as the algorithms of Tsitsiklis and Van Roy~\cite{TsiRoy99},
Longstaff and Schwarz~\cite{LoSch99} or its nonparametric version used
for \textit{optstop}. However, it should be mentioned that our new
algorithm needs much more time to compute its results: For one
of the $100$ values computed for Table~\ref{tab1}, it needs
approximately $2$ hours as opposed to $4$ minutes needed by the
\textit{optstop} algorithm.

%s5 #&#
\section{Proofs}\label{se5}

%s5.1 #&#
\subsection{\texorpdfstring{Preliminaries to the proof of Theorem \protect\ref{th1}}
{Preliminaries to the proof of Theorem 3.1}}

Once we have constructed approximations
$\hat{q}_j(z_{-\infty}^j)$
of the continuation values
$q_j(z_{-\infty}^j)$,
we can use them to construct an approximation
\[
\hat{\tau} = \inf \bigl\{ j \in\{ 0, 1, \dots, L\} \dvt g_j
\bigl(Z_1^j \bigr) \geq \hat{q}_{j}
\bigl(Z_{-\infty}^j \bigr) \bigr\}
\]
of the optimal stopping time $\tau^*$.

As our next lemma shows, the errors of the estimates $\hat{q}_j$
determine the quality of the constructed stopping time.

%le5.1 #&#
\begin{lemma}
\label{le3}
Assume $\hat{q}_L=0$. Then
\[
\EXP \bigl\{ g_{\tau^*} \bigl(Z_1^{\tau^*} \bigr) |
Z_{-\infty}^{-1} \bigr\} - \EXP \bigl\{ g_{\hat{\tau}}
\bigl(Z_1^{\hat{\tau}} \bigr) | Z_{-\infty}^{-1}
\bigr\} \leq \sum_{j=0}^{L-1} \EXP \bigl\{
\bigl\llvert \hat{q}_{j} \bigl(Z_{-\infty}^{j} \bigr) -
q_{j} \bigl(Z_{-\infty}^{j} \bigr) \bigr\rrvert |
Z_{-\infty}^{-1} \bigr\}.
\]
\end{lemma}
The assertion follows from a modification of the proof of
Proposition 21
in Belomestny~\cite{kohlerBel09}. For the sake of completeness,
a complete proof is given in the \hyperref[app]{Appendix}.

%s5.2 #&#
\subsection{\texorpdfstring{Proof of Theorem \protect\ref{th1}}{Proof of Theorem 3.1}}
Stationarity of $(Z_n)_{n \in\Z}$ implies that
%
%e5.1 #&#
\begin{equation}
\label{se2eq8} Z_{-\infty}^j \mbox{ and } Z_{-\infty}^{j+l}
\mbox{ have the same distribution for all } l \in\Z.
\end{equation}

In the sequel, we want to bound
\begin{eqnarray*}
V_0 - \hat{V}_{0,n} &=& \EXP \bigl\{ g_{\tau^*}
\bigl(Z_1^{\tau^*} \bigr) - g_{\hat{\tau}_n}
\bigl(Z_1^{\hat{\tau}_n} \bigr) \bigr\} = \EXP \bigl\{ \EXP \bigl\{
g_{\tau^*} \bigl(Z_1^{\tau^*} \bigr) - g_{\hat{\tau}_n}
\bigl(Z_1^{\hat{\tau}_n} \bigr) | Z_{-\infty}^{-1}
\bigr\} \bigr\}.
\end{eqnarray*}
By Lemma~\ref{le3}, we have
\begin{eqnarray*}
&& V_0 - \hat{V}_{0,n} \leq \sum
_{j=0}^{L-1} \EXP \bigl\{ \bigl\llvert
\hat{q}_{j,n} \bigl(Z_{-n+1}^j \bigr) -
q_{j} \bigl(Z_{-\infty}^j \bigr) \bigr\rrvert \bigr
\},
\end{eqnarray*}
so it suffices to show
%
%e5.2 #&#
\begin{equation}
\EXP \bigl\{ \bigl\llvert \hat{q}_{j,n} \bigl(Z_{-n+1}^j
\bigr) - q_{j} \bigl(Z_{-\infty}^j \bigr) \bigr
\rrvert \bigr\} \rightarrow0\qquad (n \rightarrow\infty)
\end{equation}
for $j \in\{0, 1, \dots, L-1\}$.

Using the definition of $\hat{q}_{j,n}$ as arithmetic mean
and the triangle inequality, we get
\begin{eqnarray*}
\EXP \bigl\{ \bigl\llvert \hat{q}_{j,n} \bigl(Z_{-n+1}^j
\bigr) - q_{j} \bigl(Z_{-\infty}^j \bigr) \bigr
\rrvert \bigr\} &=& \EXP \Biggl\{ \Biggl\llvert \frac{1}{n} \sum
_{l=1}^n \hat{q}_{j}^{(l)}
\bigl(Z_{-l+1}^j \bigr) - q_{j}
\bigl(Z_{-\infty}^j \bigr) \Biggr\rrvert \Biggr\}
\\
& \leq & \frac{1}{n} \sum_{l=1}^n
\EXP \bigl\{ \bigl\llvert \hat{q}_{j}^{(l)}
\bigl(Z_{-l+1}^j \bigr) - q_{j}
\bigl(Z_{-\infty}^j \bigr) \bigr\rrvert \bigr\}
\\
& \stackrel{(\ref{se2eq8})} {=} & \frac{1}{n} \sum
_{l=1}^n \EXP \bigl\{ \bigl\llvert
\hat{q}_{j}^{(l)} \bigl(Z_{1}^{j+l}
\bigr) - q_{j} \bigl(Z_{-\infty}^{j+l} \bigr) \bigr
\rrvert \bigr\}
\\
& = & \EXP \Biggl\{ \frac{1}{n} \sum_{l=1}^n
\bigl\llvert \hat{q}_{j}^{(l)} \bigl(Z_{1}^{j+l}
\bigr) - q_{j} \bigl(Z_{-\infty}^{j+l} \bigr) \bigr
\rrvert \Biggr\}.
\end{eqnarray*}

Because of the Cauchy--Schwarz inequality, it suffices to show
%
%e5.3 #&#
\begin{equation}
\EXP \Biggl\{ \frac{1}{n} \sum_{l=1}^n
\bigl\llvert \hat{q}_{j}^{(l)} \bigl(Z_{1}^{j+l}
\bigr) - q_{j} \bigl(Z_{-\infty}^{j+l} \bigr) \bigr
\rrvert^2 \Biggr\} \rightarrow0
\end{equation}
$(n \rightarrow\infty)$
for all $j \in\{0, \dots, L-1\}$.
And because of boundedness
of the estimates and of $q_j$ this in turn follows from
%
%e5.4 #&#
\begin{equation}
\label{pth1eq17} \frac{1}{n} \sum_{l=1}^n
\bigl\llvert \hat{q}_{j}^{(l)} \bigl(Z_{1}^{j+l}
\bigr) - q_{j} \bigl(Z_{-\infty}^{j+l} \bigr) \bigr
\rrvert^2 \rightarrow0
\end{equation}
in probability
for all $j \in\{0, \dots, L-1\}$.

The idea is now to use techniques from Section 27.5
(in particular Corollary 27.1) in Gy\"orfi \textit{et al.}~\cite{GyKoKrWa02}.
We have in mind definitions (\ref{eqw*}), (\ref{eqw**}) and (\ref{eqw***}),
also $\hat{q}_L^{(n)}:=0$, and
define estimates $\hat{m}_j^{(n)}(z_1^n, \cdot)$
of $m_j:=q_j$ using
realizations $z_1, \dots, z_n$ of
$Z_1$, \dots, $Z_n$, with arguments
$u_1$, \dots, $u_{n+j}$. We start
with
\[
\hat{m}_{L}^{(n)} \bigl(z_1^n ;
\cdot \bigr):=0.
\]
Given $\hat{m}^{(n)}_{j+1}(z_1^n; \cdot)$
for $j \in\{0,1,\dots,L-1\}$ we define
$\hat{m}^{(n)}_{j}(z_1^n; \cdot)$
as follows.

We start with defining $\hat{m}_{j,n,(k,h)}(z_1^n; \cdot)$
with parameters
$k \in\N$ and $h>0$
using local averaging (around $u_1^{n+j}$)
by
%
%e5.5 #&#
\begin{eqnarray}
\label{eqw****} && \hat{m}_{j,n,(k,h)} \bigl(z_1^n
;u_{1}^{n+j} \bigr)
\nonumber
\\
&&\quad := \sum_{i=k+1}^{n-j-1} \max \bigl\{
g_{j+1} \bigl(z_{i+1}^{i+j+1} \bigr),
\hat{m}_{j+1}^{(i)} \bigl(z_1^i;
z_1^{i+j+1} \bigr) \bigr\}
\\
&&\hspace*{48pt} {} \cdot { K \biggl( \frac{
u_{n-k}^{n+j}
-
z_{i-k}^{i+j}
}{h} \biggr) } \bigg/{ \sum
_{l=k+1}^{n-j-1} K \biggl( \frac{
u_{n-k}^{n+j}
-
z_{l-k}^{l+j}
}{h}
\biggr) }.
\nonumber
\end{eqnarray}
Here we set
\[
\hat{m}_{j,n,(k,h)} \bigl(z_1^n; \cdot \bigr) :=0
\]
for $k \geq n-j-1$.

For $(k,h) \in\P$ (where $\P$ is the parameter set in the definition
of the estimate), define the cumulative loss of the estimate with
parameter $(k,h)$ by
\begin{eqnarray*}
\hat{L}_{n,j}(k,h)&:=& \hat{L}_{n,j} \bigl(k,h;
z_1^{n-1}; u_1^{n+j} \bigr)
\\
&\hspace*{2.8pt}=& \frac{1}{n} \sum_{i=1}^{n-1}
\bigl( \hat{m}_{j,i,(k,h)} \bigl( z_1^i;u_{1}^{i+j}
\bigr) - \max \bigl\{ g_{j+1} \bigl(u_{i+1}^{i+j+1}
\bigr), \hat{m}_{j+1}^{(i)} \bigl(z_1^i;u_1^{i+j+1}
\bigr) \bigr\} \bigr)^2.
\end{eqnarray*}

Put $c := 8 B^2$ (where $B$ is the bound on the gain functions),
let $(p_{k,r})_{k,r}$
be the probability distribution used in the definition of the estimate
(which satisfies $p_{k,r}>0$ for all $k,r \in\N$)
and define weights, which depend on these cumulative losses, by
\[
w_{n,k,r}^{(j)}:= w_{n,k,r}^{(j)}
\bigl(z_1^{n-1}; u_1^{n+j} \bigr)
=p_{k,r} \cdot \mathrm{e}^{-n \hat{L}_{n,j}(k,h_r)/c}
\]
and their normalized values by
\[
v_{n,k,r}^{(j)}:= v_{n,k,r}^{(j)}
\bigl(z_1^{n-1}; u_1^{n+j} \bigr) =
\frac{w_{n,k,r}^{(j)}}{
\sum_{s,t=1}^\infty
w_{n,s,t}^{(j)}
}.
\]
The estimate
$
\hat{m}_{j}^{(n)}
$
is defined
as the convex combination of all estimates
$\hat{m}_{j,n,(k,h_r)}$
using weights
$v_{n,k,r}^{(j)}$, that is,
$
\hat{m}_{j}^{(n)}
$
is defined by
\[
\hat{m}_{j}^{(n)} \bigl(z_1^{n};
u_1^{n+j} \bigr) := \sum_{k,r=1}^\infty
v_{n,k,r}^{(j)} \cdot \hat{m}_{j,n,(k,h_r)}
\bigl(z_1^{n};u_{1}^{n+j} \bigr).
\]

By using a backward induction with respect to $j$ starting with
$L$, it is easy to see that we have
\[
\hat{m}_{j,(k,h)}^{(n)}=\hat{m}_{j,n,(k,h)},\qquad
\hat{q}_j^{(n)} \bigl(z_{-n+1}^j \bigr) =
\hat{m}^{(n)}_j \bigl( z_{-n+1}^0;
z_{-n+1}^j \bigr),
\]
further
\[
\hat{Q}_{n,j} (k,h) = \hat{Q}_{n,j} \bigl(z_{-n+1}^j,k,h
\bigr) = \hat{L}_{n,j} \bigl(k,h;z_{-n+1}^{-1};z_{-n+1}^j
\bigr).
\]
Thus, (\ref{pth1eq17}) means
%
%e5.6 #&#
\begin{equation}
\label{pth1eq1} \frac{1}{n} \sum_{l=1}^n
\bigl\llvert \hat{m}_{t}^{(l)} \bigl(Z_{1}^l;
Z_{1}^{l+t} \bigr) - m_{t} \bigl(Z_{-\infty}^{l+t}
\bigr) \bigr\rrvert^2 \rightarrow0
\end{equation}
in probability
for all $t \in\{0,1, \dots, L\}$, which we show
by backward induction with respect to $t$.

We start with $t=L$ in which the assertion is trivial since
\[
\hat{m}^{(l)}_L=0 \quad\mbox{and}\quad m_L=0
\]
for all $l \in\N$.

Assume now that (\ref{pth1eq1}) holds
for $t=j+1$ for some $j \in\{0,1,\dots, L-1\}$.
We have to show that in this case it is also valid
for $t=j$.

Set
\begin{eqnarray*}
L_n (\hat{m}_j) &:=& \frac{1}{n} \sum
_{l=1}^{n-1} \bigl| \hat{m}_j^{(l)}
\bigl(Z_1^l ; Z_1^{l+j} \bigr) -
\max \bigl\{ g_{j+1} \bigl(Z_{l+1}^{l+j+1} \bigr),
m_{j+1} \bigl(Z_{-\infty}^{l+j+1} \bigr) \bigr\}
\bigr|^2,
\\
L_n (\hat{m}_{j,\cdot,(k,h)}) &:=& \frac{1}{n} \sum
_{l=1}^{n-1} \bigl| \hat{m}_{j,l,(k,h)}
\bigl(Z_1^l ; Z_1^{l+j} \bigr) -
\max \bigl\{ g_{j+1} \bigl(Z_{l+1}^{l+j+1} \bigr),
m_{j+1} \bigl(Z_{-\infty}^{l+j+1} \bigr) \bigr\}
\bigr|^2,
\\
\hat{L}_n (\hat{m}_j) &:=& \frac{1}{n} \sum
_{l=1}^{n-1} \bigl| \hat{m}_j^{(l)}
\bigl(Z_1^l ; Z_1^{l+j} \bigr) -
\max \bigl\{ g_{j+1} \bigl(Z_{l+1}^{l+j+1} \bigr),
\hat{m}_{j+1}^{(l)} \bigl( Z_1^l ;
Z_1^{l+j+1} \bigr) \bigr\} \bigr|^2
\end{eqnarray*}
and
\begin{eqnarray*}
\hat{L}_n (\hat{m}_{j,\cdot,(k,h)}) &:=& \hat{L}_{n,j}
\bigl(k,h;Z_1^{n-1} ; Z_1^{n+j} \bigr)
\\
& :=& \frac{1}{n} \sum_{l=1}^{n-1} \bigl|
\hat{m}_{j,l,(k,h)} \bigl(Z_1^l ;
Z_1^{l+j} \bigr) - \max \bigl\{ g_{j+1}
\bigl(Z_{l+1}^{l+j+1} \bigr), \hat{m}_{j+1}^{(l)}
\bigl( Z_1^l ; Z_1^{l+j+1} \bigr)
\bigr\} \bigr|^2.
\end{eqnarray*}
By Lemma 27.3 in Gy\"orfi \textit{et al.}~\cite{GyKoKrWa02}, we get
%
%e5.7 #&#
\begin{equation}
\label{pth1eq4} \hat{L}_n (\hat{m}_j) \leq
\inf_{k,r \in\N} \biggl( \hat{L}_n (\hat{m}_{j,\cdot,(k,h_r)}) - c
\cdot\frac{\ln p_{k,r}}{n} \biggr).
\end{equation}
Set
\[
L_j^* := \EXP \bigl\{ \bigl| m_j \bigl(Z_{-\infty}^j
\bigr) - \max \bigl\{ g_{j+1} \bigl(Z_{1}^{j+1}
\bigr), m_{j+1} \bigl( Z_{-\infty}^{j+1} \bigr) \bigr\}
\bigr|^2 \bigr\}.
\]
In order to show (\ref{pth1eq1}), we use the following lemma which
we prove directly after the end of this proof.
%
%le5.2 #&#
\begin{lemma}
\label{le4}
Let $j \in\{0,\dots,L-1\}$. If (\ref{pth1eq1}) holds for $t=j+1$, then
%
%e5.8 #&#
\begin{equation}
\label{pth1eq2} L_n (\hat{m}_j) \rightarrow
L_j^*\qquad \mbox{in probability.}
\end{equation}
\end{lemma}

We use (\ref{pth1eq2}) to show
(\ref{pth1eq1}) for $t=j$. To do this, we proceed as in
the proof of Corollary 27.1
in Gy\"orfi \textit{et al.}~\cite{GyKoKrWa02}.
Consider the following decomposition:
\begin{eqnarray*}
&& \bigl( \hat{m}_{j}^{(l)} \bigl(Z_1^l;
Z_1^{l+j} \bigr) - \max \bigl\{ g_{j+1}
\bigl(Z_{l+1}^{l+j+1} \bigr), m_{j+1}
\bigl(Z_{-\infty}^{l+j+1} \bigr) \bigr\} \bigr)^2
\\
&& \quad= \bigl( \hat{m}_{j}^{(l)} \bigl(Z_1^l;
Z_1^{l+j} \bigr) - m_{j} \bigl(Z_{-\infty}^{l+j}
\bigr) \bigr)^2 + \bigl( m_{j} \bigl(Z_{-\infty}^{l+j}
\bigr) - \max \bigl\{ g_{j+1} \bigl(Z_{l+1}^{l+j+1}
\bigr), m_{j+1} \bigl(Z_{-\infty}^{l+j+1} \bigr) \bigr\}
\bigr)^2
\\
&&\qquad {}+ 2 \cdot \bigl( \hat{m}_{j}^{(l)}
\bigl(Z_1^l; Z_1^{l+j} \bigr) -
m_{j} \bigl(Z_{-\infty}^{l+j} \bigr) \bigr)
\\
&& \qquad \quad{}\cdot \bigl( m_{j} \bigl(Z_{-\infty}^{l+j}
\bigr) - \max \bigl\{ g_{j+1} \bigl(Z_{l+1}^{l+j+1}
\bigr), m_{j+1} \bigl(Z_{-\infty}^{l+j+1} \bigr) \bigr\}
\bigr).
\end{eqnarray*}
By (\ref{pth1eq2}), we know
\[
\frac{1}{n} \sum_{l=1}^n \bigl(
\hat{m}_{j}^{(l)} \bigl(Z_1^l;
Z_1^{l+j} \bigr) - \max \bigl\{ g_{j+1}
\bigl(Z_{l+1}^{l+j+1} \bigr), m_{j+1}
\bigl(Z_{-\infty}^{l+j+1} \bigr) \bigr\} \bigr)^2
\rightarrow L_j^*
\]
in probability.
Furthermore, by the ergodic theorem we have
\[
\frac{1}{n} \sum_{l=1}^n \bigl(
m_{j} \bigl(Z_{-\infty}^{l+j} \bigr) - \max \bigl\{
g_{j+1} \bigl(Z_{l+1}^{l+j+1} \bigr), m_{j+1}
\bigl(Z_{-\infty}^{l+j+1} \bigr) \bigr\} \bigr)^2
\rightarrow L_j^* \qquad\mbox{a.s.}
\]
Hence, it suffices to show
%
%e5.9 #&#
\begin{eqnarray}
\label{eqw*****} && \frac{1}{n} \sum_{l=1}^n
\bigl( \hat{m}_{j}^{(l)} \bigl(Z_1^l;
Z_1^{l+j} \bigr) - m_{j} \bigl(Z_{-\infty}^{l+j}
\bigr) \bigr)
\nonumber
\\
%[-8pt]\\[-8pt]
&&\hspace*{24pt} {}\cdot \bigl( m_{j} \bigl(Z_{-\infty}^{l+j}
\bigr) - \max \bigl\{ g_{j+1} \bigl(Z_{l+1}^{l+j+1}
\bigr), m_{j+1} \bigl(Z_{-\infty}^{l+j+1} \bigr) \bigr\}
\bigr)
\\
&&\quad \rightarrow0 \qquad\mbox{a.s.}
\nonumber
\end{eqnarray}
The random variables
\begin{eqnarray*}
&& \bigl( \hat{m}_{j}^{(l)} \bigl(Z_1^l;
Z_1^{l+j} \bigr) - m_{j} \bigl(Z_{-\infty}^{l+j}
\bigr) \bigr) \cdot \bigl( m_{j} \bigl(Z_{-\infty}^{l+j}
\bigr) - \max \bigl\{ g_{j+1} \bigl(Z_{l+1}^{l+j+1}
\bigr), m_{j+1} \bigl(Z_{-\infty}^{l+j+1} \bigr) \bigr\}
\bigr)
\end{eqnarray*}
are martingale differences because of $m_j=q_j$, (\ref{le2eq2}),
stationarity and dependence of the first factor on $Z_{-\infty}^{l+j}$
(not on $Z_{-\infty}^{l+j+1}$), and they are bounded by $4B^2$.
Therefore (\ref{eqw*****})
is a consequence of Theorem A.6 in Gy\"orfi \textit{et al.}~\cite{GyKoKrWa02}
(which we apply with $c_i=1$).

%s5.3 #&#
\subsection{\texorpdfstring{Proof of Lemma \protect\ref{le4}}{Proof of Lemma 5.2}}

By $|\max\{a,b\}-\max\{a,c\}| \leq|b-c|$ $(a,b,c \in\R)$
and (\ref{pth1eq1}) for $t=j+1$, we get
\begin{eqnarray*}
&& \frac{1}{n} \sum_{l=1}^{n-1} \bigl|
\max \bigl\{ g_{j+1} \bigl(Z_{l+1}^{l+j+1} \bigr),
m_{j+1} \bigl(Z_{-\infty}^{l+j+1} \bigr) \bigr\}
\\
&& \hspace*{25pt} {} - \max \bigl\{ g_{j+1} \bigl(Z_{l+1}^{l+j+1}
\bigr), \hat{m}_{j+1}^{(l)} \bigl( Z_1^l
; Z_1^{l+j+1} \bigr) \bigr\} \bigr|^2
\\
&&\quad \leq \frac{1}{n} \sum_{l=1}^{n-1}
\bigl| \hat{m}_{j+1}^{(l)} \bigl( Z_1^l ;
Z_1^{l+j+1} \bigr) - m_{j+1} \bigl(Z_{-\infty}^{l+j+1}
\bigr) \bigr|^2
\\
&&\quad\rightarrow0 \qquad\mbox{in probability.}
\end{eqnarray*}
Using
\begin{eqnarray*}
&& \frac{1}{n} \sum_{l=1}^{n-1}
|a_l-b_l|^2 - \frac{1}{n} \sum
_{l=1}^{n-1} |a_l-c_l|^2
\\
&&\quad= \frac{1}{n} \sum_{l=1}^{n-1}
(a_l-b_l+a_l-c_l)
\cdot(c_l-b_l)
\\
&& \quad\leq \Biggl( \frac{1}{n} \sum_{l=1}^{n-1}
(a_l-b_l+a_l-c_l)^2
\Biggr)^{1/2} \cdot \Biggl( \frac{1}{n} \sum
_{l=1}^{n-1} (c_l-b_l)^2
\Biggr)^{1/2}
\end{eqnarray*}
and the boundedness of the gain functions we see that this implies
%
%e5.10 #&#
\begin{equation}
\label{pth1eq*} L_n (\hat{m}_j) - \hat{L}_n
(\hat{m}_j) \rightarrow0 \quad\mbox{and} \quad L_n (
\hat{m}_{j,\cdot,(k,h)}) - \hat{L}_n (\hat{m}_{j,\cdot,(k,h)})
\rightarrow0
\end{equation}
in probability.
Hence for an arbitrary subsequence $(n_l)_l$ of $(n)_n$, we
find a subsubsequence $(n_{l_s})_s$ of $(n_l)_l$ such that
we have with probability one
%
%e5.11 #&#
\begin{eqnarray}
\label{pth1eq5} \mathop{\lim\sup}_{s \rightarrow\infty} L_{n_{l_s}} (
\hat{m}_j) &=& \mathop{\lim\sup}_{s \rightarrow\infty} \hat{L}_{n_{l_s}}
(\hat{m}_j)
\nonumber
\\
&\stackrel{(\ref{pth1eq4}) } {\leq}& \mathop{\lim\sup}_{s \rightarrow\infty}
\inf_{k,r \in\N} \biggl( \hat{L}_{n_{l_s}} (\hat{m}_{j,\cdot,(k,h_r)}) - c
\cdot\frac{\ln p_{k,r}}{n_{l_s}} \biggr)
\nonumber
\\[-8pt]
\\[-8pt]
&\leq& \inf_{k,r \in\N} \mathop{\lim\sup}_{s \rightarrow\infty} \biggl(
\hat{L}_{n_{l_s}} (\hat{m}_{j,\cdot,(k,h_r)}) - c \cdot\frac{\ln p_{k,r}}{n_{l_s}}
\biggr)
\nonumber
\\
&=& \inf_{k,r \in\N} \mathop{\lim\sup}_{s \rightarrow\infty} L_{n_{l_s}} (
\hat{m}_{j,\cdot,(k,h_r)}).
\nonumber
\end{eqnarray}
Of course, this relation also holds if we replace $(n_{l_s})_s$
by any of its subsequences (which we will do later in the proof).

Next, we analyze $L_n (\hat{m}_{j,\cdot,(k,h_r)})$.
According to (\ref{eqw****}), we have
\begin{eqnarray*}
&& \hat{m}_{j,n,(k,h)} \bigl(Z_1^n
;v_{-n+1}^j \bigr)
\\
&&\quad = { \sum_{i=k+1}^{n-j-1} \max \bigl\{
g_{j+1}\bigl(Z_{i+1}^{i+j+1}\bigr),
\hat{m}_{j+1}^{(i)}\bigl(Z_1^i;
Z_1^{i+j+1}\bigr) \bigr\} \cdot K \biggl( \frac{
v_{-k}^{j}
-
Z_{i-k}^{i+j}
}{h}
\biggr)}
\\
&&\qquad \bigg/{ \sum_{i=k+1}^{n-j-1} K \biggl(
\frac{
v_{-k}^{j}
-
Z_{i-k}^{i+j}
}{h} \biggr) }
\\
&&\quad = \frac{A_n}{C_n} + \frac{B_n-A_n}{C_n},
\end{eqnarray*}
where
\begin{eqnarray*}
A_n &:=& \frac{1}{n-j-k-1} \sum_{i=k+1}^{n-j-1}
\max \bigl\{ g_{j+1} \bigl(Z_{i+1}^{i+j+1} \bigr),
m_{j+1} \bigl(Z_{-\infty}^{i+j+1} \bigr) \bigr\} \cdot K
\biggl( \frac{
v_{-k}^{j}
-Z_{i-k}^{i+j}
}{h} \biggr),
\\
B_n &:=& \frac{1}{n-j-k-1} \sum_{i=k+1}^{n-j-1}
\max \bigl\{ g_{j+1} \bigl(Z_{i+1}^{i+j+1} \bigr),
\hat{m}_{j+1}^{(i)} \bigl(Z_1^i;
Z_1^{i+j+1} \bigr) \bigr\} \cdot K \biggl( \frac{
v_{-k}^{j}
-
Z_{i-k}^{i+j}
}{h}
\biggr)
\end{eqnarray*}
and
\[
C_n := \frac{1}{n-j-k-1} \sum_{i=k+1}^{n-j-1}
K \biggl( \frac{
v_{-k}^{j}
-
Z_{i-k}^{i+j}
}{h} \biggr).
\]
By the ergodic theorem,
we get
\[
A_n \rightarrow \EXP \biggl\{ \max \bigl\{ g_{j+1}
\bigl(Z_1^{j+1} \bigr), m_{j+1} \bigl(
Z_{-\infty}^{j+1} \bigr) \bigr\} \cdot K \biggl( \frac{v_{-k}^j-Z_{-k}^j}{h}
\biggr) \biggr\} \qquad \mbox{a.s.}
\]
and
\[
C_n \rightarrow \EXP \biggl\{ K \biggl( \frac{v_{-k}^j-Z_{-k}^j}{h} \biggr)
\biggr\} \qquad \mbox{a.s.}
\]
If we use the continuity of the kernel function, we can even apply
an ergodic theorem in the separable Banach space of continuous functions
vanishing at infinity (with supremum norm) and get that the
almost sure convergence
of $A_n$ and $C_n$ is uniformly with respect to $v_{-k}^j$
(cf., e.g., Krengel~\cite{Kr85}, Chapter 4, Theorem 2.1).

Furthermore, using the triangle inequality,
\[
\bigl|\max\{a,b\}-\max\{a,c\}\bigr| \leq|b-c| \qquad (a,b,c \in\R),
\]
and the Cauchy--Schwarz inequality we can conclude
\begin{eqnarray*}
|B_n-A_n| &\leq& \frac{1}{n-j-k-1} \sum
_{i=k+1}^{n-j-1} \bigl\llvert \hat{m}_{j+1}^{(i)}
\bigl(Z_1^{i};Z_1^{i+j+1} \bigr) -
m_{j+1} \bigl(Z_{-\infty}^{i+j+1} \bigr) \bigr\rrvert \cdot
K \biggl( \frac{
v_{-k}^{j}
-
Z_{i-k}^{i+j}
}{h} \biggr)
\\
& \leq &\sqrt{ \frac{1}{n-j-k-1} \sum_{i=k+1}^{n-j-1}
\bigl\llvert \hat{m}_{j+1}^{(i)} \bigl(Z_1^{i};Z_1^{i+j+1}
\bigr) - m_{j+1} \bigl(Z_{-\infty}^{i+j+1} \bigr) \bigr
\rrvert^2 }
\\
&&{} \cdot \sqrt{ \frac{1}{n-j-k-1} \sum_{i=k+1}^{n-j-1}
K \biggl( \frac{
v_{-k}^{j}
-
Z_{i-k}^{i+j}
}{h} \biggr)^2 }.
\end{eqnarray*}
By the ergodic theorem, the second factor on the right-hand side above
converges to
\[
\sqrt{ \EXP \biggl\{ K \biggl( \frac{v_{-k}^j-Z_{-k}^j}{h} \biggr)^2
\biggr\} } < \infty
\]
with probability one
(where we have again uniform convergence with respect to $v_{-k}^j$),
and the first factor converges
in probability
to zero by (\ref{pth1eq1}) for $t=j+1$.
Because of $K \geq c \cdot I_{S_{0,r}}$ for suitable $c>0$, $r>0$,
where $S_{0,r}$ is the ball in $\R^{j+k+1}$ centered at $0$ with
radius $r$, we have
%
%e5.12 #&#
\begin{equation}
\label{pth1eq6} \EXP \biggl\{ K \biggl( \frac{v_{-k}^j-Z_{-k}^j}{h} \biggr) \biggr\} \geq
c \cdot\PROB_{Z_{-k}^j} \bigl( v_{-k}^j +
S_{0,r \cdot h} \bigr) >0 \qquad\PROB_{Z_{-k}^j}\!\mbox{-almost
everywhere}
\end{equation}
(cf., e.g., Gy\"orfi \textit{et al.}~\cite{GyKoKrWa02}, pages 499, 500). (If $K>0$
everywhere, then (\ref{pth1eq6}) also holds everywhere.)
Therefore,
\[
\frac{B_n-A_n}{C_n} \rightarrow0 \qquad\mbox{in probability }
\PROB_{Z_{-k}^j}\!\mbox{-almost everywhere},
\]
from which we get
\[
\hat{m}_{j,n,(k,h)} \bigl(Z_1^n
;v_{-n+1}^j \bigr) \rightarrow m_{j,(k,h)}
\bigl(v_{-k}^j \bigr)\qquad \mbox{in probability}
\]
$\PROB_{Z_{-k}^j}$-almost everywhere, where
\begin{eqnarray*}
&& m_{j,(k,h)}(v_{-k},\dots, v_j)
\\
&&\quad = { \EXP \biggl\{ \max \bigl\{ g_{j+1}\bigl(Z_1^{j+1}
\bigr), m_{j+1}\bigl( Z_{-\infty}^{j+1} \bigr) \bigr\}
\cdot K \biggl( \frac{Z_{-k}^j-v_{-k}^j}{h} \biggr) \biggr\} } \bigg/{ \EXP \biggl\{ K \biggl(
\frac{Z_{-k}^j-v_{-k}^j}{h} \biggr) \biggr\} }.
\end{eqnarray*}
Let $\epsilon>0$ be arbitrary and set
\[
S_{\epsilon} = \biggl\{ v_{-k}^j \in\R^{k+j+1}
\dvt \EXP \biggl\{ K \biggl( \frac{v_{-k}^j-Z_{-k}^j}{h} \biggr) \biggr\} >\epsilon \biggr
\}.
\]
By (\ref{pth1eq6}), we know
\[
\PROB_{Z_{-k}^j} (S_{\epsilon}) \rightarrow1 \qquad(\epsilon
\rightarrow0).
\]
Since the numerators and the denominators above converge uniformly
with respect to $v_{-k}^j$ and since the limit of the denominators
is greater than $\epsilon$ on $S_{\epsilon}$, we know in addition
%
%e5.13 #&#
\begin{equation}
\label{pth1eq8} \sup_{v_{-n+1}, \dots, v{-k-1} \in\R, v_{-k}^j \in S_{\epsilon}} \bigl\llvert \hat{m}_{j,n,(k,h)}
\bigl(Z_1^n ;v_{-n+1}^j \bigr) -
m_{j,(k,h)} \bigl(v_{-k}^j \bigr) \bigr\rrvert
\rightarrow0
\end{equation}
in probability. In the sequel, we want to use this to show
%
%e5.14 #&#
\begin{eqnarray}
\label{pth1eq7} && L_{n,j}(\hat{m}_{j,\cdot,(k,h)})
\nonumber
\\
&& \quad= \frac{1}{n} \sum_{i=1}^{n-1}
\bigl( \hat{m}_{j,i,(k,h)} \bigl( Z_1^i;Z_{1}^{i+j}
\bigr) - \max \bigl\{ g_{j+1} \bigl(Z_{i+1}^{i+j+1}
\bigr), m_{j+1} \bigl(Z_{-\infty}^{i+j+1} \bigr) \bigr\}
\bigr)^2
\\
&& \quad\rightarrow \EXP \bigl\{\bigl | m_{j,(k,h)} \bigl(Z_{-k}^j
\bigr) - \max \bigl\{ g_{j+1} \bigl(Z_{1}^{j+1}
\bigr), m_{j+1} \bigl( Z_{-\infty}^{j+1} \bigr) \bigr\}
\bigr|^2 \bigr\}
\nonumber
\end{eqnarray}
in probability. To do this, we observe first that the ergodic theorem
implies
\begin{eqnarray*}
&& L_{n,j}(m_{j,(k,h)})
\\
&&\quad = \frac{1}{n} \sum_{i=1}^{n-1}
\bigl( m_{j,(k,h)} \bigl( Z_{i-k}^{i+j} \bigr) - \max
\bigl\{ g_{j+1} \bigl(Z_{i+1}^{i+j+1} \bigr),
m_{j+1} \bigl(Z_{-\infty}^{i+j+1} \bigr) \bigr\}
\bigr)^2
\\
&&\quad \rightarrow \EXP \bigl\{ \bigl| m_{j,(k,h)} \bigl(Z_{-k}^j
\bigr) - \max \bigl\{ g_{j+1} \bigl(Z_{1}^{j+1}
\bigr), m_{j+1} \bigl( Z_{-\infty}^{j+1} \bigr) \bigr\}
\bigr|^2 \bigr\}
\end{eqnarray*}
almost surely. Because of boundedness of the payoff function, we
have in addition
\begin{eqnarray*}
&& \bigl\llvert L_{n,j}(\hat{m}_{j,\cdot,(k,h)}) - L_{n,j}(m_{j,(k,h)})
\bigr\rrvert
\\
&& \quad\leq c_1 \cdot \frac{1}{n} \sum
_{i=1}^{n-1} \bigl| \hat{m}_{j,i,(k,h)} \bigl(
Z_1^i;Z_{1}^{i+j} \bigr) -
m_{j,(k,h)} \bigl( Z_{i-k}^{i+j} \bigr) \bigr|
\\
&& \quad\leq c_2 \cdot\frac{1}{n} \sum
_{i=1}^{n-1} I_{S_{\epsilon}^c} \bigl(Z_{i-k}^{i+j}
\bigr)
\\
&&\qquad{} + c_1 \cdot \frac{1}{n} \sum
_{i=1}^{n-1} \sup_{v_{-i+1}, \dots, v_{-k-1} \in\R, v_{-k}^j \in S_{\epsilon}} \bigl\llvert
\hat{m}_{j,i,(k,h)} \bigl(Z_1^i ;v_{-i+1}^j
\bigr) - m_{j,(k,h)} \bigl(v_{-k}^j \bigr) \bigr
\rrvert
\\
&& \quad\rightarrow c_2 \cdot\PROB_{Z_{-k}^j}
\bigl(S_{\epsilon}^c \bigr)
\end{eqnarray*}
in probability
by (\ref{pth1eq8}) and by the ergodic theorem. By letting
$\epsilon\rightarrow0$, we get (\ref{pth1eq7}).
And by replacing $(n_{l_s})$ by a suitable subsequence of
$(n_{l_s})_s$, we can assume
w.l.o.g. even that (\ref{pth1eq7}) holds for almost sure
convergence if we replace $n$ by $n_{l_s}$ in (\ref{pth1eq7}).

Next, we use Lemma 24.8 in Gy\"orfi \textit{et al.}~\cite{GyKoKrWa02} which implies
\[
m_{j,(k,h)} \bigl(z_{-k}^j \bigr) \rightarrow
m_{j,k} \bigl(z_{-k}^j \bigr)\qquad
\PROB_{Z_{-k}^j}\!\mbox{-almost everywhere}
\]
for $h \rightarrow0$,
where
\[
m_{j,k} \bigl(z_{-k}^j \bigr) := \EXP \bigl\{
\max \bigl\{ g_{j+1} \bigl(Z_1^{j+1} \bigr),
m_{j+1} \bigl(Z_{-\infty}^{j+1} \bigr) \bigr\}
|Z_{-k}^j=z_{-k}^j \bigr\}.
\]
And by the martingale convergence theorem,
we have
\[
m_{j,k} \bigl(Z_{-k}^j \bigr) \rightarrow
m_j \bigl(Z_{-\infty}^j \bigr)\qquad \mbox{a.s.}
\]
for $k \rightarrow\infty$
(since the almost sure limit $X$ of the left-hand side satisfies
\[
\int_A X \,\mathrm{d}P = \int_A
m_j \bigl(Z_{-\infty}^j \bigr) \,\mathrm{d}P
\]
for all $A \in\F(Z_{-k}^j)$ and all $k \in\N$, cf., e.g.,
Chapter 32.4A in Lo\`eve~\cite{L77} for more general results in this respect).
From this, we conclude by dominated convergence
\begin{eqnarray*}
&& \mathop{\lim\sup}_{s \rightarrow\infty} L_{n_{l_s}}(\hat{m}_j)
\\
&&\quad\stackrel{(\ref{pth1eq5})} {\leq} \inf_{k,r \in\N} \mathop{\lim
\sup}_{s \rightarrow\infty} L_{n_{l_s}} (\hat{m}_{j,\cdot,(k,h_r)})
\\
&&\quad\stackrel{(\ref{pth1eq7})} {=} \inf_{k,r \in\N} \EXP \bigl\{ \bigl|
m_{j,(k,h_r)} \bigl(Z_{-k}^j \bigr) - \max \bigl\{
g_{j+1} \bigl(Z_1^{j+1} \bigr), m_{j+1}
\bigl( Z_{-\infty}^{j+1} \bigr) \bigr\} \bigr|^2 \bigr\}
\\
&&\quad\leq L_j^*\qquad \mbox{a.s.}
\end{eqnarray*}
Because of
\[
\mathop{\lim\inf}_{n \rightarrow\infty} L_n(\hat{m}_j) \geq
L_j^*\qquad \mbox{a.s.}
\]
(cf., e.g., Section 27.5 in Gy\"orfi \textit{et al.}~\cite{GyKoKrWa02})
this completes the proof of (\ref{pth1eq2}).

%%apA #&#
%%\begin{appendix}\label{app}
%%
%%sB #&#
%%\section*{Proof of Lemma \protect\ref{le3}}

\begin{appendix}\label{app}
\setcounter{equation}{0}

\section*{\texorpdfstring{Appendix: Proof of Lemma \protect\ref{le3}}{Appendix: Proof of Lemma 5.1}}

Set
\[
\hat{\tau}_t^* = \inf \bigl\{ s \geq t+1 \dvt \hat{q}_s
\bigl(Z_{-\infty}^s \bigr) \leq g_s
\bigl(Z_1^s \bigr) \bigr\}
\]
and let $\F_t$ be the $\sigma$-algebra generated by
$Z_{-\infty}^t$. In the sequel, we prove
%
%eB.1 #&#
\begin{eqnarray}
\label{ple3eq1} && \EXP \bigl\{ g_{\tau_{t-1}^*} \bigl(Z_1^{\tau_{t-1}^*}
\bigr) - g_{\hat{\tau}^*_{t-1}} \bigl(Z_1^{\hat{\tau}^*_{t-1}} \bigr) |
\F_{t-1} \bigr\} %\nonumber\\
%&&
\leq \sum
_{k=t}^{L-1} \EXP \bigl\{ \bigl| \hat{q}_k
\bigl(Z_{-\infty}^k \bigr) - q_k
\bigl(Z_{-\infty}^k \bigr) \bigr| | \F_{t-1} \bigr\}\quad
\end{eqnarray}
for $t \in\{0, \dots, L\}$, from which we get
the assertion of Lemma~\ref{le3} by setting $t=0$.

We prove (\ref{ple3eq1}) by induction. The assertion
is trivial for $t=L$ (since
$\tau_{L-1}^*=L=\hat{\tau}^*_{L-1}$). Assume that
(\ref{ple3eq1}) holds for $t \in\{ s+1, \dots, L \}$
for some $s \in\{0,1, \dots, L-1\}$.
In the sequel we prove that in this case it also holds for $t=s$.
To do this, we use
\begin{eqnarray*}
&& \EXP \bigl\{ g_{\tau_{t-1}^*} \bigl(Z_1^{\tau_{t-1}^*} \bigr) -
g_{\hat{\tau}^*_{t-1}} \bigl(Z_1^{\hat{\tau}^*_{t-1}} \bigr) |
\F_{t-1} \bigr\}
\\
&& \quad= \sum_{k=t}^{L-1} \EXP \bigl\{
\bigl(g_{\tau_{t-1}^*} \bigl(Z_1^{\tau_{t-1}^*} \bigr) -
g_{\hat{\tau}^*_{t-1}} \bigl(Z_1^{\hat{\tau}^*_{t-1}} \bigr) \bigr) \cdot
1_{
\{
\hat{\tau}^*_{t-1}=k, \tau_{t-1}^*>k
\}
} | \F_{t-1} \bigr\}
\\
&&\qquad {}+ \sum_{k=t}^{L-1} \EXP \bigl\{
\bigl( g_{\tau_{t-1}^*} \bigl(Z_1^{\tau_{t-1}^*} \bigr) -
g_{\hat{\tau}^*_{t-1}} \bigl(Z_1^{\hat{\tau}^*_{t-1}} \bigr) \bigr) \cdot
1_{
\{
\hat{\tau}_{t-1}^*>k, \tau_{t-1}^*=k
\}
} | \F_{t-1} \bigr\}
\\
&&\quad = \sum_{k=t}^{L-1} \EXP \bigl\{
\bigl( g_{\tau_{k}^*} \bigl(Z_1^{\tau_{k}^*} \bigr) -
g_{k} \bigl(Z_1^{k} \bigr) \bigr) \cdot
1_{
\{
\hat{\tau}_{t-1}^*=k, \tau_{t-1}^*>k
\}
} | \F_{t-1} \bigr\}
\\
&& \qquad{}+ \sum_{k=t}^{L-1} \EXP \bigl\{
\bigl( g_{k} \bigl(Z_1^{k} \bigr) -
q_k \bigl(Z_{-\infty}^k \bigr) \bigr) \cdot
1_{
\{
\hat{\tau}_{t-1}^*>k, \tau_{t-1}^*=k
\}
} | \F_{t-1} \bigr\}
\\
&&\qquad {}+ \sum_{k=t}^{L-1} \EXP \bigl\{
\bigl( q_k \bigl(Z_{-\infty}^k \bigr) -
g_{\hat{\tau}_{k}^*} \bigl(Z_1^{\hat{\tau}_{k}^*} \bigr) \bigr) \cdot
1_{
\{
\hat{\tau}_{t-1}^*>k, \tau_{t-1}^*=k
\}
} | \F_{t-1} \bigr\}
\\
&& \quad= T_1 + T_2 + T_3,
\end{eqnarray*}
where we have used that
$\hat{\tau}_{t-1}^*=\hat{\tau}_{k}^*$ on
$\{\hat{\tau}_{t-1}^*>k\}$ and that
$\tau_{t-1}^*={\tau}_{k}^*$ on
$\{{\tau}_{t-1}^*>k\}$.
The random variables
\[
1_{
\{
\hat{\tau}_{t-1}^*=k, \tau_{t-1}^*>k
\}
} \quad\mbox{and} \quad 1_{
\{
\hat{\tau}_{t-1}^*>k, \tau_{t-1}^*=k
\}
}
\]
are $\F_k$-measurable, hence we get by Lemma~\ref{le2}
\begin{eqnarray*}
T_1 &=& \sum_{k=t}^{L-1} \EXP
\bigl\{ \bigl( \EXP \bigl\{ g_{\tau_{k}^*} \bigl(Z_1^{\tau_{k}^*}
\bigr) | \F_k \bigr\} - g_{k} \bigl(Z_1^{k}
\bigr) \bigr) \cdot 1_{
\{
\hat{\tau}_{t-1}^*=k, \tau_{t-1}^*>k
\}
} | \F_{t-1} \bigr\}
\\
&=& \sum_{k=t}^{L-1} \EXP \bigl\{ \bigl(
q_k \bigl(Z_{-\infty}^k \bigr) - g_{k}
\bigl(Z_1^{k} \bigr) \bigr) \cdot 1_{
\{
\hat{\tau}_{t-1}^*=k, \tau_{t-1}^*>k
\}
} |
\F_{t-1} \bigr\}
\\
& \leq& \sum_{k=t}^{L-1} \EXP \bigl\{
\bigl( q_k \bigl(Z_{-\infty}^k \bigr) -
\hat{q}_k \bigl(Z_{-\infty}^k \bigr) \bigr) \cdot
1_{
\{
\hat{\tau}_{t-1}^*=k, \tau_{t-1}^*>k
\}
} | \F_{t-1} \bigr\}, %\label{ple3eq2}
\end{eqnarray*}
since $\hat{\tau}_{t-1}^*=k$ implies
\[
g_{k} \bigl(Z_1^{k} \bigr) \geq
\hat{q}_k \bigl(Z_{-\infty}^k \bigr).
\]
Similarly, $\hat{\tau}_{t-1}^*>k$ implies
\[
g_{k} \bigl(Z_1^{k} \bigr) <
\hat{q}_k \bigl(Z_{-\infty}^k \bigr),
\]
from which we can conclude
\[
T_2 \leq \sum_{k=t}^{L-1} \EXP
\bigl\{ \bigl( \hat{q}_k \bigl(Z_{-\infty}^k \bigr) -
q_k \bigl(Z_{-\infty}^k \bigr) \bigr) \cdot
1_{
\{
\hat{\tau}_{t-1}^*>k, \tau_{t-1}^*=k
\}
} | \F_{t-1} \bigr\}.
\]
Finally, we have by Lemma~\ref{le2}
\begin{eqnarray*}
T_3 &=& \sum_{k=t}^{L-1} \EXP
\bigl\{ \EXP \bigl\{ q_k \bigl(Z_{-\infty}^k \bigr) -
g_{\hat{\tau}_{k}^*} \bigl(Z_1^{\hat{\tau}_{k}^*} \bigr) | \F_k
\bigr\} \cdot 1_{
\{
\hat{\tau}_{t-1}^*>k, \tau_{t-1}^*=k
\}
} | \F_{t-1} \bigr\}
\\
&=& \sum_{k=t}^{L-1} \EXP \bigl\{ \EXP
\bigl\{ g_{\tau_{k}^*} \bigl(Z_1^{\tau_{k}^*} \bigr) -
g_{\hat{\tau}_{k}^*} \bigl(Z_1^{\hat{\tau}^*_{k}} \bigr) | \F_k
\bigr\} \cdot 1_{
\{
\hat{\tau}_{t-1}^*>k, \tau_{t-1}^*=k
\}
} | \F_{t-1} \bigr\},
\end{eqnarray*}
and by using the induction hypothesis we get
\begin{eqnarray*}
T_3 &\leq& \sum_{k=t}^{L-1}
\EXP \Biggl\{ \sum_{j=k+1}^{L-1} \EXP \bigl\{ \bigl|
\hat{q}_j \bigl(Z_{-\infty}^j \bigr) -
q_j \bigl(Z_{-\infty}^j \bigr)\bigr| | \F_k
\bigr\} \cdot 1_{
\{
\hat{\tau}_{t-1}^*>k, \tau_{t-1}^*=k
\}
} | \F_{t-1} \Biggr\}
\\
&=& \sum_{k=t}^{L-1} \EXP \Biggl\{ \sum
_{j=k+1}^{L-1} \bigl| \hat{q}_j
\bigl(Z_{-\infty}^j \bigr) - q_j
\bigl(Z_{-\infty}^j \bigr)\bigr| \cdot 1_{
\{
\hat{\tau}_{t-1}^*>k, \tau_{t-1}^*=k
\}
} |
\F_{t-1} \Biggr\}
\\
&=& \sum_{j=t+1}^{L-1} \EXP \Biggl\{ \bigl|
\hat{q}_j \bigl(Z_{-\infty}^j \bigr) -
q_j \bigl(Z_{-\infty}^j \bigr)\bigr| \cdot \sum
_{k=t}^{j-1} 1_{
\{
\hat{\tau}_{t-1}^*>k, \tau_{t-1}^*=k
\}
} | \F_{t-1}
\Biggr\}
\\
&=& \sum_{k=t+1}^{L-1} \EXP \Biggl\{\bigl |
\hat{q}_k \bigl(Z_{-\infty}^k \bigr) -
q_k \bigl(Z_{-\infty}^k \bigr)\bigr| \cdot \sum
_{j=t}^{k-1} 1_{
\{
\hat{\tau}_{t-1}^*>j, \tau_{t-1}^*=j
\}
} | \F_{t-1}
\Biggr\}.
\end{eqnarray*}
Summarizing the above results, we get the assertion.

\end{appendix}

\section*{Acknowledgement}

The authors thank anonymous referees for various suggestions
substantially improving the presentation of the results.

%suskaldyti doi

% imsref loaded by aiste.veprauskaite, 2012-08-03 11:24:47

\printhistory

\end{document}